\documentclass[12pt]{amsart}  
\usepackage{amssymb,amsmath,amsthm,amscd}
\usepackage[all]{xy}

\numberwithin{equation}{section}
\newtheoremstyle{fancy1}{10pt}{10pt}{\itshape}{12pt}{\textsc\bgroup}{.\egroup}{8pt}{
}
\newtheoremstyle{fancy2}{10pt}{10pt}{}{12pt}{\itshape}{.}{8pt}{ }

\theoremstyle{fancy1}

\newtheorem{cor}[equation]{Corollary}

\newtheorem{prop}[equation]{Proposition}
\newtheorem*{prop*}{Proposition}
\newtheorem{thm}[equation]{Theorem}

\newtheorem*{prob*}{Problem}
\newtheorem{main}{Theorem}
\newtheorem*{mainc}{Corollary C}
\newtheorem*{main*}{Theorem}

\newtheorem*{cor*}{Corollary}

\newtheorem*{problem*}{Problem}

\theoremstyle{fancy2}

\newtheorem{rem}[equation]{Remark}
\newtheorem*{rem*}{Remark}

\newcommand{\cref}[1]{Corollary~\ref{#1}}

\newcommand{\pref}[1]{Proposition~\ref{#1}}

\newcommand{\tref}[1]{Theorem~\ref{#1}}

\newcommand{\e}{\epsilon}

\newcommand{\CP}{\mathbb{C\mkern1mu P}}

\newcommand{\C}{{\mathbb{C}} }

\newcommand{\Z}{{\mathbb{Z}} }

\newcommand{\Ra}{{\mathbb{Q}} }
\newcommand{\Sph}{\mathbb{S}}

\newcommand{\G}{\ensuremath{\operatorname{G}} }
\newcommand{\D}{\ensuremath{\operatorname{\D}} }
\newcommand{\SO}{\ensuremath{\operatorname{SO}} }

\newcommand{\Sp}{\ensuremath{\operatorname{Sp}} }
\newcommand{\U}{\ensuremath{\operatorname{U}} }
\newcommand{\SU}{\ensuremath{\operatorname{SU}} }

\renewcommand{\S}{\ensuremath{\operatorname{S}} }

\newcommand{\fg}{{\mathfrak{g}} }
\newcommand{\fk}{{\mathfrak{k}} }
\newcommand{\fh}{{\mathfrak{h}} }

\def\con#1=#2(#3){#1 \equiv #2 \bmod{#3}}

\newcommand{\diag}{\ensuremath{\operatorname{diag}} }

\DeclareMathOperator{\Id}{Id}

\newcommand{\no}{\noindent}

\newcommand{\B}{\mathcal{B}}

\begin{document}
\title{On the topology of positively curved Bazaikin spaces}  

\author{Luis A. Florit}
\address{IMPA: Est. Dona Castorina 110, 22460-320, Rio de Janeiro,
Brazil}
\email{luis@impa.br}
\author{Wolfgang Ziller}
\address{University of Pennsylvania: Philadelphia, PA 19104, USA}
\email{wziller@math.upenn.edu}
\thanks{The first author was supported by CNPq-Brazil
 and the second author  by the Francis
J. Carey Term Chair,  by a grant from the National Science
Foundation, and by the Clay Institute.}

\maketitle

Manifolds with positive sectional curvature have been of interest
since the beginning of global Riemannian geometry, as illustrated by
the theorem of Bonnet-Myers, the Berger-Klingenberg sphere theorem,
and the Synge theorem. Surprisingly, for compact simply connected
manifolds, no obstructions are known to distinguish between manifolds
(that admit a metric) with non-negative curvature and those of
positive curvature. On the other hand, among known examples, the class
of non-negatively curved manifolds is much larger. They include, among
others, any homogeneous space and more generally the so called
biquotients, i.e quotients of a compact Lie groups $G$ by a subgroup
of $G\times G$ acting on $G$ from the left and right simultaneously.
A biinvariant metric on $G$ induces a non-negatively curved metric on
the quotient.

The difficulty of finding new examples of manifolds with positive
curvature may be illustrated by the fact that all known examples are
constructed in this way, allowing more generally left invariant
metrics on $G$. But only very few manifolds of this type have positive
curvature and, apart from the rank one symmetric spaces, exist only in
low dimensions. They consist of certain homogeneous spaces in
dimensions $6,7,12,13$ and $24$ due to Berger \cite{Be}, Wallach
\cite{Wa}, and Aloff-Wallach \cite{AW}, and of biquotients in
dimensions $6,7$ and $13$ due to Eschenburg \cite{E1},\cite{E2} and
Bazaikin \cite{Ba}.

Aside from 4
isolated manifolds, this list divides into two infinite
families. The first in dimension $7$ consisting of the Eschenburg
biquotients
$$
{\mathcal E}=\diag(z^{k_1}, z^{k_2},
z^{k_3})\backslash \SU(3)/ \diag(z^{l_1}, z^{l_2}, z^{l_3})^{-1},
$$
with {\small$\sum$} $\!k_i=$ {\small$\sum$} $\!l_i$, which include
the infinite sub-family of homogeneous Aloff-Wallach spaces
${\mathcal W}=\SU(3)/\diag(z^{l_1},z^{l_2},z^{l_3})$,
{\small$\sum$} $\!l_i=0$.
The second one exists in dimension 13 and consists of the Bazaikin
biquotients
$$
\B_q =
\diag(z^{q_0},z^{q_1},z^{q_2},z^{q_3},z^{q_4},z^{q_5})\backslash
\SU(6)/\Sp(3),
$$
where $q=(q_0,\dots,q_5)$ are odd integers with {\small$\sum$}
$\!q_i=0$. When the restrictions on $q$ that are necessary for this
biquotient to be a manifold are satisfied, we call $\B_q$ a
{\it Bazaikin space}. If, in addition, the condition for a certain
metric on $\B_q$ to have positive sectional curvature holds, we call
$\B_q$ a {\it positively curved Bazaikin space}; see Section 1.
\smallskip

Topological properties of the Eschenburg spaces $\mathcal E$ have
been studied extensively. A complete classification of their
homeomorphism and diffeomorphism types has been obtained in
\cite{KS} and \cite{K1},\cite{K2}. Using these, it was shown that
occasionally positively curved homogeneous Aloff-Wallach spaces
(\cite{KS}), and more frequently positively curved Eschenburg
biquotients (\cite{CEZ}), can be homeomorphic to each other, but not
diffeomorphic. Furthermore, Eschenburg spaces are frequently
homotopy equivalent and occasionally homeomorphic and even
diffeomorphic to homogeneous Aloff-Wallach spaces without being
Aloff-Wallach spaces by definition; cf. \cite{CEZ}. See also
\cite{AMP1}, \cite{AMP2}, and \cite{Sh}.

\vskip .2cm
In this paper we study topological properties of the Bazaikin spaces
and we will see that their behaviour is quite different from the
Eschenburg spaces. A further difference is that there is only one
Bazaikin space which is homogeneous by definition, the Berger space
$\SU(5)/\Sp(2)\S^1$ (see \cite{Be}), which corresponds to
$q=(-5,1,1,1,1,1)$ (up to order). On the other hand, there exists
a large class of homogeneous spaces in dimension $13$. However, we
will show:

\begin{main}\label{mainB}
A Bazaikin space $\B_q$ with $q\ne (-5,1,1,1,1,1)$ is not homeomorphic
to a homogeneous space. If, in addition, $\B_q$ is positively curved,
then it is not even homotopy equivalent to a homogeneous space.
\end{main}

We point out that there are positively curved Bazaikin spaces which
have the same cohomology ring as certain homogeneous spaces, and
there are Bazaikin spaces which match all of the  homotopy invariants
(as below) of a Berger space.

In contrast to Theorem \ref{mainB}, there are positively curved
Eschenburg spaces that are homotopically equivalent but not
homeomorphic to Aloff-Wallach spaces (\cite{Sh}).

\vskip .2cm

The proof involves, among other things, the computation of certain
topological invariants. So far no classification of the Bazaikin
spaces up to homeomorphism or diffeomorphism has been obtained. In
fact, only the cohomology ring was known (\cite{Ba}). The property of
the cohomology ring that can distinguish Bazaikin spaces among
each other is the order $s$ of the finite torsion groups
$H^6(\B_q,\Z)=H^8(\B_q,\Z)=\Z_s$. Further invariants are given by the
Pontryagin classes $p_1\in H^4(\B_q,\Z)=\Z$ and
$p_2\in H^8(\B_q,\Z)=\Z_s$, and the linking form $lk\in \Z_s$.

\begin{main}\label{main2}
For a Bazaikin space $\B_q$ the order $s$ is odd and:
$$
8s=|\sigma_3| \; ,   \; lk = \pm 32\sigma_5^{-1}\in \Z_s,
$$
$$
p_1=-\sigma_2=\tfrac 1 2 \Vert q\Vert^2 \; ,  \;
8p_2 = 3p_1^2-\sigma_4 \in \Z_s,
$$
where $\sigma_i=\sigma_i(q)$ stands for the elementary
symmetric polynomial of degree $i$ in $q$.
\end{main}

Recall that $s$ and $lk$, as well as $p_1\!\!\mod 24$, are homotopy
invariants, while rational Pontryagin classes are homeomorphism
invariants. For a Bazaikin space, $p_1$ is hence a homeomorphism invariant
while $p_2$ is a diffeomorphism invariant.

\begin{mainc}
The following finiteness results hold:
\begin{itemize}
\item[$i)$] There are only finitely many positively curved Bazaikin
spaces for a given cohomology ring;
\item[$ii)$] There are only finitely many Bazaikin spaces in each
homeomorphism type.
\end{itemize}
\end{mainc}

Part $(i)$ is also true for Eschenburg spaces (\cite{CEZ}) and should
be viewed in the context of the Klingenberg-Sakai conjecture. It
states that there are only finitely many positively curved manifolds
in a given homotopy type and raises the question weather this might be
true even for a given cohomology ring.

It is not difficult to find Bazaikin spaces that match all four
invariants in Theorem~\ref{main2}, especially for small values of
$s$. However, by means of a computer program, we observed that:
\vskip 0.3cm
\centerline{\it All positively curved Bazaikin spaces with
$s\leq 10^9$,}
\centerline{\it altogether 2.130.601.485 manifolds, are
homeomorphically distinct.}
\vskip 0.3cm
\no This strongly suggests that indeed all positively
curved Bazaikin spaces are homeomorphically, or at least
diffeomorphically, distinct.

\vskip .2cm
In Section 1 we discuss preliminaries which are needed, in Section 2
we compute the homotopy invariants, while in Section 3 we compute the
Pontryagin classes and present the proof of Theorem~\ref{mainB}.
We end the paper with some computer experiment results in Section~4.

\vskip .2cm
We would like to thank C. Escher and T. Chinburg for helpful
discussions. This work was done while the second author was visiting
IMPA and he would like to thank the Institute for its hospitality.
\vskip .8cm

\section{Preliminaries}  

A biquotient can be defined in several ways. First, consider two
subgroups of $G$ defined by monomorphisms $f_1\colon H\to G$ and
$f_2\colon K\to G$. The group $H\times K$ acts on $G$ via
$(h,k)\cdot g = f_1(h)gf_2(k)^{-1}$. If this action is free, the
quotient $H\backslash G/K$ is called a {\it biquotient}. In the
literature one often finds the seemingly more general definition where
$H\times K$ is allowed to be replaced by a subgroup of $L$ of
$G\times G$. This can be reduced to the above case by rewriting the
biquotient as $\triangle G\backslash G\times G/L$. The action is free
iff $(h,k)=(e,e)$ whenever $h\in H$ and $k\in K$ are conjugate to each
other. This clearly can be further reduced to the case where we assume
that $h$ and $k$ lie in a maximal torus.

In the following we also allow the action of $H\times K$ to have
a finite ineffective kernel $\Gamma$ (necessarily embedded diagonally
in $H\times K$) and that $(H\times K)/\Gamma$ acts freely on $G$. This
again is not more general since it can be reduced to the above by
replacing $G$ with $G/\Gamma$. But the advantage is that the
topological computations become simpler when $G$ is simply connected.

\smallskip

Following \cite{EKS}, we define a Bazaikin space as
\begin{equation}\label{Ba1}
\B_q = \diag(z^{q_0},z^{q_1},z^{q_2},z^{q_3},z^{q_4},z^{q_5})
\backslash\SU(6)/\Sp(3),
\end{equation}
where $z\in \S^1=\{w\in\C:|w|=1\}$, $q=(q_0,\cdots ,q_5)$ is an
ordered set of integers with $\sum q_i=0$. We choose the embedding
$\Sp(3)\subset \SU(6)$ such that a maximal torus of $\Sp(3)$ is given
by $\diag(u,\bar{u},v,\bar{v},w,\bar{w})\subset\SU(6)$. This
biquotient action of $\S^1\times \Sp(3)$ on $\SU(6)$ has ineffective
kernel $\Gamma=\Z_2= \langle(-1, -\Id)\rangle$ and one easily sees
that the action of $(\S^1\times \Sp(3))/\Gamma$ is free if and only if
\begin{equation}\label{e:gcd}
\text{\it all $q_i$'s are odd\ \ and\ \ }\gcd(q_{\tau (1)}+q_{\tau
(2)},q_{\tau (3)}+ q_{\tau (4)})=2,
\end{equation}
for all permutations $\tau\in S_6.$ Notice that under a reordering of
the integers $q_i$, or replacing $q$ with $-q$, one obtains the same
manifold.

\smallskip

The manifolds $\B_q$ can be written in an alternative way. Since the
symmetric space $\SU(6)/\Sp(3)$ can be described as the set of unitary
almost complex structures, one easily sees that the subgroup
$\SU(5)\subset\SU(6)$ acts transitively with isotropy group $\Sp(2)$.
Indeed, this follows from the usual normal form for such complex
structures since one can fix the first vector arbitrarily. If $\SU(5)$
is the subgroup that fixes $e_0=(1,0,\cdots,0)\in \C^6$, it follows
that $\B_q$ becomes a quotient of $\SU(5)$ under the action of
$\S^1\times \Sp(2)$:
\begin{equation}\label{Ba2}
\B_q = \diag (z^{q_1} , \dots , z^{q_5} )
\backslash \SU(5)/ \diag (z^{-q_0},A)^{-1},
\end{equation}
where $A\in \Sp(2)\subset \SU(4)\subset\SU(5)$ with the embedding of
$\Sp(2)$ and $\Sp(3)$ as described above. In this language it is
natural to describe the space by the 5 integers
$\bar{q}=(q_1,\dots,q_5)$. Accordingly, we will refer to $\B_q$ in
either way, namely, by $q\in \Z^6$ or $\bar q\in \Z^5$, and no
confusion will arise because of the different number of integers
involved. But it is important to observe that, for the $\bar q \in
\Z^5$ notation, aside from ordering and sign, replacing one of the
integers by the negative of their sum also gives the same space. E.g.,
the Berger space can be written as $(1,1,1,1,1)$ or as $(-5,1,1,1,1)$.

This action again has an ineffective kernel
\mbox{$\Z_2=\langle(-1,-\Id)\rangle$}. To make this second description
into a true biquotient, we use the diffeomorphism
$\SU(5)\simeq \U(5)/\diag (w,{\rm Id})$ to rewrite it as
$$
\B_{\bar{q}} = \diag (z^{q_1} , \dots , z^{q_5} )
\backslash \U(5)\,/ \diag (w,A)^{-1},
$$
with $z,w\in \S^1$ and $A\in\Sp(2)$.

In \cite{Ba} one finds the somewhat different description
\begin{equation}\label{Ba3}
\diag (z^{q_1'} , \dots , z^{q_5'} )
\backslash \U(5)\,/ \diag(1,wA)^{-1},
\end{equation}
and one sees, by writing both as quotients of $\SU(5)$,
that they are related by
\begin{equation}\label{trans}
4q_i'=q_i+q_0 \;,\; i=1,\dots ,5.
\end{equation}
Notice that one needs $ q_i\equiv \epsilon$ mod $4$ for all
$i=1,\dots,5$, for a fixed $\epsilon=\pm1$, in order to rewrite
a Bazaikin space as defined in \eqref{Ba2} into one as defined in
\eqref{Ba3}. This turns out to be a strong condition for the
positively curved Bazaikin spaces (see Section 4).

\smallskip

There are two natural fibrations associated to a Bazaikin space.
First, $\B_q$ is, by definition, the base space of a principal
circle bundle:
\begin{equation}\label{fibS1}
\S^1 \rightarrow \SU(6)/ \Sp(3) \rightarrow \B_q.
\end{equation}
But notice that the $\S^1$ action is only free modulo the ineffective
$\Z_2$ kernel. Second, the total space of this fibration has a further
natural fibration
\begin{equation}\label{fibM}
\Sph^5\rightarrow\SU(6)/\Sp(3) \rightarrow \Sph^9,
\end{equation}
To see this fibration, we use the fact that
$\SU(6)/\Sp(3)=\SU(5)/\Sp(2)$. One then obtains a fibration from
the inclusions $\Sp(2)\subset \SU(4)\subset \SU(5)$ and the fact
that $\SU(4)/\Sp(2)=\SO(6)/\SO(5)=\Sph^5$ and
$\SU(5)/\SU(4)=\Sph^9$.

\smallskip

In order to obtain metrics with positive curvature, the symmetric metric
needs to be modified. Each Bazaikin space has six natural biquotient
metrics associated to it. First, one can rewrite a space as in
\eqref{Ba1} in six different ways as in \eqref{Ba2} by choosing
different base points $e_0$ along the coordinate axes. This
corresponds to removing one of the six integers in $q$. For a space
as in \eqref{Ba2} we have the natural left invariant metrics on
$\SU(5)$ obtained from the biinvariant metric by scaling with scale
$<1$ in the direction of $\U(4)\subset\SU(5)$ embedded in the last 4
coordinates. Since $\S^1\times\Sp(2)$ acts by isometries in this
modified metric, it induces a submersion metric on $\B_q$.

In \cite{Zi} it was shown (see \cite{EKS} for a published version)
that the necessary and sufficient conditions for one of these
metrics to have positive sectional curvature is that
\begin{equation}\label{pos}
\text{\it there exists } 0\leq i_0 \leq 5 \text{ \it such that
$q_i+q_j>0$ {\rm (}or $<0${\rm )}},\; \text{\it for all }
i_0\ne i < j \ne i_0.
\end{equation}
It is easy to check that at most one of these 6 metrics has positive
curvature. In the original paper \cite{Ba}, there are four
conditions ensuring positive sectional curvature for a  metric as
defined in \eqref{Ba3}. But notice that condition (c) in
\cite[p.1069]{Ba} implies the other three and that in the
translation \eqref{trans} condition (c) corresponds to $q_i>0 $ for
$i=1,\dots,5$. Bazaikin's condition is again restrictive. In fact,
around $7\%$ of all examples satisfy both of Bazaikin's extra
assumptions; see Section 4.

\smallskip

We point out that the change \eqref{Ba1} in the description of the
Bazaikin spaces turns out to be more than cosmetic. It simplifies the
expressions for the invariants and the proof of our theorems
significantly. Notice that for the symmetric polynomials one has
$\sigma_i(q)=\sigma_i(\bar{q})-\sigma_1(\bar{q})\sigma_{i-1}(\bar{q})$.

\section{Homotopy Invariants}  

In this section we discuss the topology of the Bazaikin spaces. We
start with the cohomology ring. Although already computed in
\cite{Ba}, we present a proof here since it is much simpler and some
of the information obtained will be used again in the computation of
the Pontryagin classes.

\begin{prop}[\cite{Ba}]\label{coh}
The Bazaikin space $\B_q$ is simply connected and the non-vanish\-ing
cohomology groups are given by
\begin{eqnarray*}
H^0(\B_q)=\, H^2(\B_q)=\, H^4(\B_q)=H^9(\B_q)=
 H^{11}(\B_q)=\,  H^{13}(\B_q)=\, \Z, \,\\ \text{ and } \,
H^6(\B_q)=H^8(\B_q) = \Z_s,\hspace{100pt}
\end{eqnarray*}
with $8s=|\sigma_3(q)|$. The ring structure is determined by the fact
that if $u\in H^2(\B_q)$ is a generator, then $u^i$ is a generator
of $H^{2i}(\B_q)$ for $i=2,3,4$.
\end{prop}

\begin{proof}
As discussed in Section 1, let $G=\SU(6)$, $H=\S^1$ and  $K=\Sp(3)$
with embeddings $f_1\colon H\to G$ given by
$f_1(z)=\diag (z^{q_0},\dots, z^{q_5})\in \SU(6)$ and
$f_2\colon K\to G$ the canonical embedding of $\Sp(3)\subset\SU(6)$.
Furthermore, assume that the biquotient action of $(H\times K)/\Gamma$
on $G$ is free, where $\Gamma=\Z_2=\langle(-1,-\Id)\rangle$.

For brevity, $LM$, in the following, will always denote the direct
product $L\times M$ of the connected Lie groups $L$ and $M$.
Furthermore, $B_L$ will denote the classifying space of $L$. Choose
a contractible space $E$ on which $GG$ acts freely. Hence $HK$ and $H$
or $K$ also act freely with quotient their respective classifying
space.

The fibration in \eqref{fibS1} induces the fibration
\begin{equation}\label{BS1}
G/K \rightarrow \B_q\rightarrow B_{H'},
\end{equation}
where $H'=H/\{\pm 1\}$ is the free action of $\S^1$ on $G/K$. On the
other hand, the Serre spectral sequence of the fibration \eqref{fibM}
implies that \mbox{$H^*(G/K,\Z)=H^*(\Sph^5\times \Sph^9,\Z)$}. Hence
in the spectral sequence of the fibration \eqref{BS1} one has only two
non-vanishing differentials:
$$
d_{6}\colon E_5^{0,5}=\Z\to E_5^{6,0}=\Z \; , \; d_{10}\colon
E_{11}^{0,9}=\Z\to E_{11}^{10,0}=\Z_s.
$$
The differential $d_6$ is multiplication by some integer $s$ and,
together with Poincare duality,  it follows that the cohomology ring
of $\B_q$ is as claimed in the Theorem with
$H^6(\B_q)=H^8(\B_q)=\Z_s$. It remains to determine the integer $s$.
For this purpose we use the naturality of differentials in the
following commutative diagram of fibrations:

\begin{picture}(100,160)\label{dddd1}  
\put(139,130){$G\times E$}
\put(250,130){$G/K$}
\put(190,133){\vector(1,0){48} }
\put(155,117){\vector(0,-1){20} }
\put(262,117){\vector(0,-1){20} }
\put(130,80){$ G \times_{HK} E$ }
\put(257,80){$\B_q$}
\put(190,83){\vector(1,0){48} }
\put(207,72){$\pi$}
\put(155,65){\vector(0,-1){20} }
\put(160,55){{$(\varphi_{H},\varphi_{K})$ } }
\put(262,67){\vector(0,-1){20} }
\put(270,55){{$B_\Delta$ } }
\put(130,30){$B_{H}\times B_{K}$}
\put(255,30){$B_{H'}$}
\put(190,33){\vector(1,0){48} }
\put(190,3){\it \small Diagram 1.}
\end{picture}  
\vskip 0.5cm
\noindent
Here $HK$ acts freely on $G\times E$ via a diagonal action consisting
of the left and right action on $G$ and the free action on $E$. The
projection onto the second factor $G\times_{HK}E\to B_{H}\times B_{K}$
is the classifying map of this $HK$ principle bundle and defines the
left hand side fibration. All horizontal maps are the natural
projections.

We first claim that the differential $d_6^{\; '}$ in the left hand
side fibration determines $d_6$. To~see this, observe that the
projection $G \times E\to G/K$ induces an isomorphism in
dimension~$5$. This follows using the edge homomorphism in the
spectral sequence of \mbox{$G\!\times\!E\!\to\!G/K\!\to\!B_K$} since
the differentials on $E_5^{0,5}=\Z$ vanish for dimension reasons.
Furthermore, since \mbox{$H\to H'$} is given by $z\to z^2$, the
homomorphism $H^*(B_{H'})\to H^*(B_{H})$ is, via transgression,
multiplication by $2$ in dimension~$2$ and hence by $2^i$ in
dimension $2i$. This describes the induced map
$H^*(B_{H'})\to H^*(B_{H}\times B_{K})$ and shows that $d_6^{\; '}$
determines $d_6$.

To compute the differential $d_6^{\; '}$, we use naturality with
respect to another commutative diagram of fibrations, modifying the
method used in \cite{E3} for Eschenburg spaces:

\begin{picture}(100,160)\label{dddd2}  
\put(133,130){$G\times E$}
\put(245,130){$G\times E$}
\put(180,133){\vector(1,0){55} }
\put(145,117){ \vector(0,-1){20} }
\put(258,117){  \vector(0,-1){20} }
\put(123,80){$ G \times_{HK} E$ }
\put(250,80){$B_{\Delta\G}=G\times_{GG}E$}
\put(180,83){\vector(1,0){55} }
\put(200,72){$\varphi_{G}$}
\put(145,65){  \vector(0,-1){20} }
\put(150,55){ {$(\varphi_{H},\varphi_{K})$ } }
\put(258,65){  \vector(0,-1){20} }
\put(265,55){ {$B_\Delta$ } }
\put(124,30){$B_{H}\times B_{K}$}
\put(240,30){$B_{G}\times B_{G}$}
\put(180,33){\vector(1,0){55} }
\put(190,22){$B_{(f_1,f_2)}$}
\put(180,-3){\it \small Diagram 2.}
\end{picture}  
\vskip 0.5cm
\noindent
The right hand side fibration comes from the $GG$ principal bundle
associated to the free diagonal action of left and right
multiplication on $G$ and the free action on $E$. Notice that
$G\times_{GG}E=E/G=B_G$.

Since all Lie groups involved have no torsion in their cohomology, all
remaining computations can be done with integer coefficients, and all
cohomology groups are to be understood with $\Z$ coefficients. It is
well known that $H^*(G)=H^*(\SU(6))$ is the exterior algebra
$\Lambda(y_2,\cdots , y_6)$ and $H^*(B_{G})$ the polynomial algebra
$P[\bar{y}_2,\cdots , \bar{y}_6]$. Here, $y_i$ has degree $2i-1$ and
its transgression $\bar{y}_i$ has degree $2i$. If $T_{G}$ is the
maximal torus of $G$ with coordinates
$(t_1,\dots , t_6) \; , \; \sum t_i=0$ we identify $t_i$, by abuse of
notation, with the elements $t_i\in H^1(T_{G})$ and hence
$\bar{t}_i\in H^2(B_{T_{G}})$. We then have
$H^*(B_{T_{G}})=P[\bar{t}_1\dots , \bar{t}_6]$ and
$H^*(B_{G})= H^*(B_{T_{G}})^{W_{G}}$, where $W_{G}$ is the Weyl group
of $G$. A basis of the algebra of $W_{G}$\,--\,invariant elements is
given by the elementary symmetric polynomials
$\sigma_i(\bar{t})=\sigma_i(\bar{t}_1,\dots , \bar{t}_6)$ and hence we
can choose $\bar{y}_i=\sigma_i(\bar{t})$.

Similarly, for $T_{K}$ we use the coordinates $(s_1,s_2,s_3)$
and hence $H^*(B_{K})=H^*(B_{T_{K}})^{W_{K}}=P[\bar{s}_1,
\bar{s}_2,\bar{s}_3]^{W_{K}}=P[b_1,b_2,b_3]$, with
$b_i=\sigma_i (\bar{s}_1^2, \bar{s}_2^2,\bar{s}_3^2)$. Finally,
$H^*(B_{H})=P[\bar{u}]$, where $u$ is the coordinate of the
circle $H$.

Under the embeddings $f_1\colon T_{H}\to T_{G}$ with
$f_1(u)=(q_0u,\dots, q_5u)$ we have \mbox{$f_1^*(t_i)=q_i u$}
and hence
$B_{f_1}^*(\bar{y}_i)=
B_{f_1}^*(\sigma_i(\bar{t}))=\sigma_i(q)\bar{u}^i$.
Furthermore, $f_2\colon T_{K}\to T_{G}$ is
given by $f_2(s_1,s_2,s_3)=(s_1,-s_1,s_2,-s_2,s_3,-s_3)$ and
clearly $B_{f_2}^*(\bar{y}_{2i+1})=0$ since $H^{2i+1}(B_K)=0$.

As was shown in \cite{E3}, the differentials in the spectral sequence
for the fibration $B_{\Delta G}\to B_{GG}$ are given by
$$
d_{2i}(y_i)=\bar{y}_i\otimes 1 - 1\otimes \bar{y}_i,
$$
where we identify $y_i\in H^*(G)$ with  elements in $ E^{0,*}_2$ and
$H^*(B_{G}\times B_{G})$  with elements in $E^{*,0}_2$. By naturality,
the differentials in the spectral sequence for the left hand side
fibration in Diagram 2 are given by
$d_{2i}^{\; '}(y_i)=B_{f_1}^*(\bar{y}_i)-B_{f_2}^*(\bar{y}_i)$
and hence $d_6^{\;'}(y_3)= \sigma_3(q)\bar{u}^3$. We conclude that
$d_6^{\;'}\colon E_5^{0,5}=\Z\to E_5^{6,0}=\Z$ is multiplication by
$\sigma_3(q)$.

This proves our claim since, as we observed earlier, $H^*(B_{H'})\to
H^*(B_{H})$ is multiplication by $8$ in dimension $6$.
\end{proof}

\smallskip

In \cite{CEZ} it was shown that for a given positively curved
Eschenburg space, there are only finitely many other  positively
curved Eschenburg spaces with the same cohomology ring.  As claimed
in Corollary C,  the analogous result for Bazaikin spaces is true as
well:
\medskip

{\it Proof of Corollary C (i)}.\ \
The only variable cohomology group is $H^6(\B_q)=H^8(\B_q)=\Z_{s}$.
We reorder the integers and change the sign of $q$ so that
$q_0\leq q_1 \leq \dots\leq q_5$ and $|q_2|\leq q_3$. The positive
curvature condition is then equivalent to $q_1+q_2> 0$. We can rewrite
$s$ as
\begin{align}\label{order}
8s = -\sigma_3(q) = &\ (q_1+q_2)^2(q_3+q_4+q_5)+(q_1+q_2)
((q_3+q_4+q_5)^2+q_1q_2)\\
&+(q_3+q_4)(q_4+q_5)(q_3+q_5),\nonumber
\end{align}
which clearly implies the desired claim.
\qed
\bigskip

A more subtle homotopy invariant is given by the linking form which is
defined as follows. Let $X$ be a simply connected manifold whose
cohomology ring agrees with the cohomology ring of a Bazaikin space
for some integer $s$. Throughout this note, we will call such a space
a {\it homological Bazaikin space}. Consider the Bockstein homomorphism
$\beta \colon H^5(X,\Z_s)=\Z_s\to H^6(X,\Z)=\Z_s$, associated to the
short exact sequence $0\to\Z\overset{\cdot s}{\to}\Z\to\Z_s\to 0$.
The long exact sequence
$$
\cdots \to H^5(X,\Z)\to H^5(X,\Z_s)\overset{\beta}{\longrightarrow}
H^6(X,\Z)\to\cdots
$$
implies that $\beta$ is an isomorphism since $H^5(X,\Z)=0$.
The linking form is then given by
$$
L\colon H^6(X,\Z)\times H^8(X,\Z) \to \Z_s \; , \;
L(a,b)=(\beta^{-1}(a)\cup b)([X]),
$$
for a given choice of orientation class $[X]$. $L$ is clearly
determined up to sign by $lk:=L(u^3,u^4)\in \Z_s$, where $u$ is a
generator of $H^2(X,\Z)$. We now claim that

\begin{thm}\label{linking}
The linking form of $\B_q$ is given by
$$
lk= \pm32\sigma_5^{-1} \in \Z_s,
$$
with $\sigma_5$ and $s$ as in Theorem B.
\end{thm}

\begin{proof}
The number $lk$ can be described in another fashion. Let $X$ be
a homological Bazaikin space. Define an $\S^1$ principal bundle $P$
over $X$ by requiring that its Euler class $e(P)$ is a generator of
$H^2(X)=\Z$. A change in sign for $e(P)$ correspond to a change of
orientation of the circle bundle. Hence $P$ is well defined up to
orientation and is determined by the homotopy type of $X$.
 From the Gysin sequence for
$\S^1\to P\to X$ and Poincare duality for $P$ it follows that
$H^*(P,\Z)=H^*(\Sph^5\times \Sph^9,\Z)$. We can now consider the
spectral sequence of the bundle
$$
P\to X\to B_{\S^1},
$$
whose isomorphism type is a homotopy invariant of $X$. As in the proof
of \pref{coh}, it is determined by only two differentials which can be
non-zero, namely,
$$
d_{6}\colon E_5^{0,5}=\Z\to E_5^{6,0}=\Z \; , \;  d_{10}\colon
E_{11}^{0,9}=\Z\to E_{11}^{10,0}=\Z_s.
$$
Here $d_6$ must be multiplication by $s$ since $H^6(X)=\Z_s$ and
$d_{10}$ must be multiplication by some $v$ with $\gcd(s,v)=1$ since
$H^{10}(X)=0$. This homotopy invariant $v$ is related to the linking
form. In \cite{Si} one finds a discussion of the corresponding
question for the Eschenburg spaces where $H^*(P,\Z)=H^*(\Sph^3\times
\Sph^5,\Z)$ and the non-zero differentials are $d_3$ and $d_5$. It is
shown there that $lk=\pm v^{-1}$. One easily sees that the same proof
carries over in our situation. It thus remains to compute the
invariant $v$.

In the case of $X=\B_q$ the total space $P$ of the circle bundle can
be chosen to be $P=\SU(6)/\Sp(3)$. Indeed, by the long homotopy
sequence of $\Sp(3)\to \SU(6)\to P$ it follows that
$\pi_1(P)=\pi_2(P)=0$ and hence $H^2(P)=0$. The Gysin sequence of the
circle bundle then implies that the Euler class is a generator of
$H^2(\B_q)$. The fibration $P\to X\to B_{\S^1}$ now agrees with the
one in \eqref{BS1}. We can hence compute the differentials in the
spectral sequence of this fibration using again Diagram 1 and
Diagram 2. As in \pref{coh}, one shows that $d_{10}$ is determined by
$d_{10}^{\;'}$ and that $H^{10}(B_{H'}) \to H^{10}(B_{H})$ is
multiplication by $32$. Furthermore, one has
$B_{f_1}^*(\bar{y}_5)= \sigma_5(q)\bar{u}^5$ and
$B_{f_2}^*(\bar{y}_5)=0$ and hence
$d_{10}^{\;'}(y_5)=\sigma_5(q)\bar{u}^5$.
As in \pref{coh}, this proves our claim.
\end{proof}

\bigskip

There is another homotopy invariant associated to a homological
Bazaikin space. To describe it, consider, as in the proof of
\tref{linking}, the total space $P$ of the circle bundle with Euler
class a generator in $H^2(X)$. As explained there, $P$ is simply
connected with $H^*(P,\Z)=H^*(\Sph^5\times \Sph^9,\Z)$. Hence one can
describe $P$ as a CW--complex
$$
P = \Sph^5 \cup e^{9}\cup \cdots,
$$
with cells in dimension $5,9,\cdots$. The attaching map $\partial
(e^{9})=\Sph^{8}\to \Sph^5$ is an element
$\alpha\in \pi_8(\Sph^5)=\Z_{24}$ (cf. \cite{Ha, To}) and hence
$\pi_8 (P)=\Z_{28}/\langle\alpha\rangle$. Since
$\pi_8 (X)=\pi_8 (P)$, the group $\pi_8 (X)$ must be a quotient of
$\Z_{24}$ as well. Thus, homological Bazaikin spaces naturally fall
into 8 different homotopy types.

For the Bazaikin spaces we claim that they belong to the class where
$\pi_8(X)=0$. Indeed, the long homotopy sequence for
$\Sp(3)\to \SU(6)\to P$:
$$
\cdots\to\pi_8(\SU(6))\to \pi_8(P)\to \pi_7(\Sp(3))
\to \pi_7(\SU(6))\to\cdots,
$$
together with the fact that $\pi_7(\SU(6))=\pi_7(\Sp(3))=\Z$ and
$\pi_8(\SU(6))=0$ (\cite{Mii}), imply that $\pi_8(P)=0$ or $\Z$ and
thus vanishes.

\smallskip

There is a second natural family of homological Bazaikin spaces. They
depend on integers $a=(a_1,\dots ,a_3)$ and $b=(b_1,\dots ,b_5)$, and
are described as a quotient \begin{equation}\label{S1Quotient}
M_{a,b}=\Sph^5\times \Sph^9/\S^1_{a,b} \end{equation} under the circle
action
$$
\big((z_1,z_2,z_3),(w_1,\dots,w_5)\big)\in\Sph^5\times\Sph^9\to
\big((e^{ia_1\theta}z_1,\dots ,e^{ia_3\theta} z_3),
(e^{ib_1\theta}w_1,\dots,e^{ib_5\theta} w_5)\big).
$$
 From the fibration $\S^1\to\Sph^5\times \Sph^9\to M_{a,b}$ it
follows that $\pi_8(M_{a,b})=\Z_{24}$. Summarizing, we have
\begin{prop}\label{pi8}
The Bazaikin spaces belong to the homotopy type $\pi_8(X)=0$ and the
manifolds $M_{a,b}$ to the homotopy type $\pi_8(X)=\Z_{24}$.
\end{prop}

In the case where $a_1=a_2=a_3$ and $b_1=\dots=b_5$ the manifold
$M_{a,b}$ is homogeneous under the transitive action of
$\SU(3)\times\SU(5)$ on $\Sph^5\times \Sph^9$ since it commutes with
the circle action.

\begin{rem}
By combining the long homotopy sequences of the above fibrations,
one easily shows that the homotopy groups $\pi_i=\pi_i(\B_q)$,
$i\leq 13$, are given by
$$
\pi_i=0 \text{ for } i=1,3,4,8,10,11, \ \pi_6=\pi_7=\pi_{13}=\Z_2,\
\pi_{12}=\Z_{360}\ \text{ and } \pi_2=\pi_5=\pi_9=\Z
$$
(cf. \cite{To}).
In particular, one sees that the Bazaikin spaces cannot be
distinguished from each other by their low dimensional homotopy groups.
\end{rem}

As we will see in the next section, there is another homotopy
invariant given by \mbox{$p_1\!\! \mod 24$} which divides the Bazaikin
spaces into two different homotopy types.

\bigskip

\section{Diffeomorphism invariants}  

The next natural set of invariants are the Pontryagin classes, which
depend on the  diffeomorphism type. We now show:
\begin{thm}\label{pont}
The Pontryagin classes $p_1$ and $p_2$ of $\B_q$ are given by
$$
p_1=- \sigma_2\in H^4(\B_q)=\Z,
$$
$$
8p_2=3\sigma_2^2-\sigma_4\in H^8(\B_q)=\Z_{s},
$$
\vspace{4pt}
with $\sigma_i$ and $s$ as in Theorem B.
\end{thm}

\begin{proof}
To compute the Pontryagin classes, we modify \cite[Theorem 4.2]{Si} to
the situation presented in the proof of \pref{coh}. As described
there, we have a projection $\pi\colon X= G \times_{HK} E\to \B_q$ and
an $H K$ principal bundle with classifying map
$(\varphi_{H},\varphi_{K})\colon X \to B_{H}\times B_{K}$.
In \cite{Si} the following vector bundles were introduced:
$\alpha_H = (G/K) \times _H \fh$ where $H$ acts on $G/K$ on the left
and on $\fh$ via the adjoint representation. Similarly
$\alpha_K = (H\backslash G) \times _K \fk$ and $\alpha_G
= (H\backslash G)\times (G/K) \times _G \fg$ where $G$ acts on
$(H\backslash G)\times (G/K)$ via
$g\cdot(Hg_1,g_2K)=(Hg_1g^{-1},gg_2K)$ and on $\fg$ via the adjoint
representation. The projection onto the first factor shows that they
are all vector bundles over $\B_q=H\backslash G/K$. Although the
action of $HK$ on $G$ is only almost free, one easily sees that the
proof of \cite[Theorem 3.2]{Si} carries over to show that the tangent
bundle $\tau$ of $\B_q$ satisfies
$$
\tau\oplus \alpha_H\oplus\alpha_K = \alpha_G.
$$

Since the action is only almost free, \cite[Theorem 4.2]{Si} cannot
be applied directly, but can be modified as follows. Define
$\bar{\alpha}_H = (G\times_K E) \times _H \fh$ with $H$ acting on
$G$ on the left and similarly $\bar{\alpha}_K = (G\times_H E)
\times _K \fk$ and
$\bar{\alpha}_G = [(G\times_H E)\times (G\times_K E)] \times _G \fg$
as vector bundles over $X$. Clearly $\pi^*(\bar{\alpha}_H)=\alpha_H$
and similarly for the others so that we obtain
$$
\pi^*(\tau)\oplus\bar{\alpha}_H\oplus\bar{\alpha}_K=\bar{\alpha}_G.
$$
We can now use the usual formula for the Pontryagin classes of a
homogeneous vector bundle $\alpha_{L}=P\times_{L}V$ associated to the
$L$ principle bundle $P\to P/L=B$:
$$
p({\alpha}_L)=1+p_1+p_2+\cdots = \varphi_{L}^*(a) \, , \,
a=\prod (1+\alpha_i^2),
$$
where $\alpha_i$ runs through the positive weights of the
representation of $L$ on $V$ and $\varphi_{L}\colon B\to B_L$ is the
classifying map of the $L$ principal bundle. Here we have identified
$\alpha_i\in H^1(T_{L})\cong H^2(B_{T_{L}})$ and hence
$a\in H^*(B_{T_{L}})^{W_{L}}\cong H^*(B_{L})$.

The vector bundles
$\bar{\alpha}_H \, , \, \bar{\alpha}_K \, , \, \bar{\alpha}_G$
are associated to a principal bundle and the weights are the roots of
the corresponding Lie group. Thus $p(\bar{\alpha}_H)=1$, and we obtain
$$
p(\pi^*(\tau))=\varphi_{G}^*(a)
\varphi_{K}^*(b^{-1}) \, , \,
a=\prod_{\alpha_i\in\Delta_{G}^+}(1+\alpha_i^2) \, , \,
b=\prod_{\beta_i\in\Delta_{K}^+}(1+\beta_i^2).
$$
The positive roots $\Delta_G^+$ for $G=\SU(6)$ are
$(t_i-t_j)\, , \, 1\le i<j\le 6$, with $\sum t_i=0$,
and a computation shows that
$$
a=1 - 12\overline y_2 + 60\overline y_2^2 + 12\overline y_4 + \cdots,
$$
where $\bar{y}_i=\sigma_i(\bar{t}_1,\cdots , \bar{t}_6)$ are the
elementary symmetric polynomials in $\bar{t}_i\in H^2(B_{T_{G}})$.
Similarly, the positive roots for $K= \Sp(3)$ are
$s_i\pm s_j \, ,\, 2s_i$ with $1\leq i < j\leq 3$ and hence
$$
b = 1 + 8b_1 + 22b_1^2 + 14b_2 + \cdots \text{\; and \; } b^{-1} = 1
- 8b_1 + 42b_1^2 - 14b_2 + \cdots,
$$
with $b_i\in H^2(B_{T_{K_1}}) $ defined as in the proof of \pref{coh}.

 From the spectral sequence computation in \pref{coh} it follows that
$\varphi_{G}^*(\bar{y}_i)=B_{f_1}^*(\bar{y}_i)=\sigma_i(q)\bar{u}^i$.
Furthermore, one shows that $B_{f_2}^*(\bar{y}_{2})=-b_1$ and
$B_{f_2}^*(\bar{y}_{4})=b_2$ and hence
$\varphi_{K}^*(b_1)=-\sigma_2\bar{u}^2$ and
$\varphi_{K}^*(b_2)=\sigma_4\bar{u}^2$. Altogether,
\begin{align*}
\varphi_{G}^*(a) \varphi_{K}^*(b^{-1})
&= (1\!-\!12\sigma_2 \bar{u}^2\!+\!12(5\sigma_2^2\!+\!\sigma_4)
\bar{u}^4\!+\!\cdots)(1\!+\!8\sigma_2\bar{u}^2\!+\!
14(3\sigma_2^2\!-\!\sigma_4)\bar{u}^4\!+\!\cdots)\\
&= 1-4 \sigma_2\bar{u}^2 +2(3\sigma_2^2-\sigma_4)\bar{u}^4 +\cdots,
\end{align*}
\no which computes the Pontryagin class for $\pi^*(\tau)$. Since
$\pi^*(p(\tau))=p(\pi^*(\tau))$, the Pontryagin classes for $\tau$
will now be determined by the homomorphisms
$\pi^*\colon H^i(\B)\to H^i(X), i=4,8$. Recall that in the proof of
\pref{coh} we showed that $\pi^*\colon H^2(B_{H'})\to H^2(B_H)$ is
multiplication by $2$ in dimension $2$. On the other hand, using
Diagram 2 we showed that $H^2(X) $ is generated by
$\bar{y}\in H^2(B_H)$, whereas Diagram 1 implies that
$H^2(\B_q)=H^2(B_{H'})$. Hence $\pi^*$ is multiplication by $4$ in
dimension $4$ and by $16$ in dimension $8$. Combining all of these
proves our claim.
\end{proof}

Notice that the first Pontryagin class, interpreted as integer, has
a well defined sign since $u^2$ is a uniquely defined generator in
$H^4(\B_q)=\Z$. Since the rational Pontryagin classes are not only
diffeomorphism invariants, but also homeomorphism invariants, so is
the integer $p_1=-\sigma_2(q)$. Furthermore, a theorem of Hirzebruch
(\cite{Hi}) says that $p_1\!\mod 24$ is a homotopy invariant. For
this, we have the following.

\begin{cor}\label{c:mod}
For any Bazaikin space $\B_q$, the following holds:
\begin{itemize}
\item[$a)$] The first Pontryagin class satisfies
$p_1 = 7\ {\rm or}\, 15\!\mod 24$. Moreover, up to order and sign,
we have:
\begin{itemize}
\item[$a.1)$] $p_1 = 15\!\mod 24$ if and only if
$q\!\!\mod 3=(1,1,1,1,1,1)$ or $(1,1,1,0,0,0);$
\item[$a.2)$] $p_1 = 7\!\mod 24$ if and only if
$q\!\!\mod 3=(1,1,1,1,-1,0);$
\end{itemize}
\item[$b)$] The order $s$ of $H^6(\B_q)=\Z_s$ satisfies
$s=\pm 1\! \mod 6.$
\end{itemize}
\end{cor}

\begin{proof}
Since all $q_i$'s are odd, by (\ref{e:gcd}) we have
$q\!\!\mod 4=\pm(1,1,1,1,1,-1)$. Then, we easily get that
$p_1=7\!\mod 8$ and that $s$ is odd. Again by (\ref{e:gcd}), the only
possibilities for $q\!\mod 3$ are $(1,1,1,1,1,1)$, $(1,1,1,0,0,0)$ or
$(1,1,1,1,-1,0)$ up to ordering and sign of~$q$. In all cases we
obtain $s=\pm1\!\mod 3$. An easy computation verifies the remaining
claims.
\end{proof}

Thus Bazaikin spaces  naturally fall into two homotopy types. It is
also known that \mbox{$p_2+2p_1^2$} is a homotopy invariant mod $5$
(\cite{Wu}). But this does not give any further information. Indeed,
using (\ref{e:gcd}) we see that $s\!\!\mod 5 =0$ if and only if
\mbox{$q\!\!\mod 5=\pm(\e,\e,\e,\e,\e,0)$}, for $\e=1,2$. This easily
gives $p_1\!\!\mod 5=p_2\!\!\mod 5 = 0$.
\bigskip

With the expression for the above topological invariants, we are in
position to give the
\vskip .3cm
{\it Proof of Theorem \ref{mainB}}.
Let $M=G/H$ be a homogeneous space which we assume is homotopy
equivalent to $\B_q$ for some $q$. Since it is simply connected, it is
well known that there exists a semisimple subgroup of $G$ which also
acts transitively. We can hence assume that $G$ is semisimple. We can
furthermore assume that $G$ is simply connected by making the action
on $G/H$ ineffective if necessary. Since $M$ is simply connected, this
implies that $H$ is connected. From the long homotopy sequence of the
fibrations \eqref{fibS1} and \eqref{fibM}, together with the fact that
the rational homotopy groups of an odd $n$-dimensional sphere are only
non-zero in dimension $n$, it follows that $\pi_i(\B_q)\otimes\Ra=0$
for $i\ne 2,5,9$ and that it is equal to $\Ra$ for $i= 2,5,9$. Hence
the same holds for $G/H$. Now recall that for every compact Lie group
$L$ one has $H^*(L,\Ra)=H^*(S^{n_1}\times \cdots \times S^{n_k},\Ra)$
where $L$ has rank $k$ and the integers $n_i$ are called the exponents
of $L$. These exponents are well known for each compact simple Lie
group (cf. \cite[Table 1.4]{KZ} for a complete list). In particular
they are all odd, $1$~occurs with multiplicity~$r$ if and only if the
center of $L$ is $r$-dimensional, and $3$ occurs with multiplicity $m$
if and only if $L$ has $m$ simple factors. Thus, from the long exact
homotopy sequence for $H\to G\to M$ tensored with $\Ra$, and the fact
that $G$ contains no $\S^1$ factors, it follows that $H=H'\cdot \S^1$
with $H'$ semisimple, and that $G$ and $H'$ have the same number of
simple factors. Furthermore, there must be an extra $\Sph^5$ and
$\Sph^9$ in the sphere decomposition of $G$ compared to the one for
$H'$. We can now use \cite[Lemma 1.4]{KZ}, which states that
a homomorphism of Lie groups $f\colon A\to B$ such that
$f_*\colon \pi_*(A)\otimes\Ra\to \pi_*(B)\otimes\Ra$ is onto must
itself be onto as a Lie group homomorphism as well. Hence, if $G$ has
at least three simple factors,
$G=G_1\times G_2 \times G_3\times\cdots$ and the extra 5 and 9-spheres
are contained in $G_1\times G_2$, the projection of $H$ onto
$G_3\times\cdots$ must be onto which implies that $G_1\times G_2$
already acts transitively on $G/H$. Hence we can assume that $G$ and
$H'$ both have one or two simple factors.

We now look at the list of exponents for simple Lie groups, and use
the usual abbreviations
$A_n=\SU(n+1)\, ,\, C_n=\Sp(n)\, ,\, D_n=\SO(2n)$. By abuse of
notation $C_n$ will also stand for $B_n=\SO(2n+1)$ since they have the
same exponents. Due to the low dimensional isomorphisms $D_3=A_3$ and
$C_1=A_1=D_1$, we can assume that $n\ge 1$ for $C_n$ , $n\ge 2 $ for
$A_n$ and $n\ge 4$ for $D_n$. Now observe that the only simple Lie
groups which contain a 5-sphere are $A_n$ for some $n\ge 2$, and the
only ones which contain a 9-sphere are $A_n$, $n\ge 4$, or $D_5$. One
easily obtains the following list of pairs $(G,H')$ which have all of
the above properties:
\begin{itemize}
\item[(I)]
$(G,H')$ is one of (a) $(A_4,C_2) $ or (b) $ (A_5,C_3)$;
\item[(II)] $(G,H')=(G_1G_2,H_1H_2)$ is one of the following:\\
$(a)\; (A_3D_5,C_2C_4)\, ,\,(b)\;  (A_2D_5,C_1C_4)\, ,\,
(c)\; (A_3A_4,C_2A_3), \text{\ \ or\ \ } (d)\;(A_2A_4,C_1A_3);$
\item[(III)]
$(G,H')=(G_1L,LH_1)$, where $(G_1,H_1)$ is one of the cases in I.
\end{itemize}

 Using the possible embeddings of classical Lie groups into each
other, as well as the low dimensional isomorphisms and the
replacement of $C_n$ by $B_n$ discussed above, we see that none of
the simple factors of $H'$  can be embedded diagonally in $G_1G_2$.
Hence, in case of II and III, $G/H'$ is automatically a product
homogeneous space $(G_1/H_1)\times (G_2/H_2)$  with the order
indicated in the list. Since we also have a circle inside $H$, we
need the further property that the normalizer of $H_1$ in $G_1$ or
the normalizer of $H_2$ in $G_2$, or both, are not finite. This
excludes Ib  and  IIa. The case Ia uniquely determines the Berger
space since $\SO(5)\subset \SU(5)$ has a finite normalizer and the
normalizer of $\Sp(2)\subset \SU(5)$ is one dimensional. This also
excludes case III, since one needs $H_1\subset L\subset G_1$ in this
case, which is only possible when $L$ contains a circle, contradicting
our assumption that $G$ is semisimple.

For the remaining three cases in II, we are left with the following
possibilities:
\begin{align*}
\tag{IIb} &G/H\, =(\SU(3)/\SU(2)\S^1)\times (\SO(10)/\SO(9))=
\CP^2\times \Sph^9,\\
\tag{IIc} &G/H\, =(\SO(6)/\SO(5))\times (\SU(5)/\SU(4)\S^1) =
\Sph^5 \times \CP^4,\\
\tag{IId} &G/H'=(\SU(3)/\SU(2))\times ( \SU(5)/\SU(4))=
\Sph^5\times \Sph^9, \text{ or }\\
\tag{IId$'$} &G/H\, =(\SU(3)/\SO(3))\times
(\SU(5)/\SU(4)\S^1)=(\SU(3)/\SO(3))\times \CP^4.
\end{align*}
The two cases where the manifold contains a $\CP^4$ factor contradict
the rational cohomology ring structure of the Bazaikin space. In IId
the circle in $H$ can be embedded diagonally at slope $(a,b)$ into
the two dimensional normalizer of $H'$ in $G$ and we can assume that
$\gcd(a,b)=1$. The case IIb can now be interpreted as the special
subcase $(a,b)=(1,0)$. Notice also that $a=0$ again contradicts the
rational cohomology ring of $\B_q$ and we can hence assume $a\ge 1$.

Altogether, we are left with the Berger space and the manifolds
$G/H=\Sph^5\times\Sph^9/S^1_{a,b}$, where the circle $S^1_{a,b}$ acts
on $\Sph^5\times\Sph^9$ as a Hopf action on each sphere, but with
relative speed $a$ on the first factor and $b$ on the second factor.
These are special cases of the manifolds $M_{a,b}$ discussed in
\pref{pi8} where we showed that the Bazaikin spaces and the manifolds
$\Sph^5\times\Sph^9/S^1_{a,b}$ belong to different homotopy types.

So far, we have shown that if a Bazaikin space has the homotopy type
of a homogeneous space, the homogeneous space can only be the Berger
space. We now need to show that the only Bazaikin space with this
property is the Berger space itself.

For a positively curved Bazaikin space, the equation \eqref{order}
easily implies that among the positively curved manifolds $\B_q$ the
value $s=5$ is only obtained for the Berger space
$\bar{q}=(1,1,1,1,1)$. This proves the first part.

For a general Bazaikin space $\B_q$, we use the homeomorphism
invariant $p_1=\tfrac 1 2 \|q\|^2$. For the Berger space we have
$p_1=15$ which means that if $\B_q$ is homeomorphic to it,
$\max_i(|q_i|)\le 5$. One now easily sees that $p_1=15$ and $s=5$ is
only possible when $\bar{q}=(1,1,1,1,1)$. This finishes the proof of
the second assertion.
\qed
\smallskip

\begin{rem*}\label{r:hom}
$a)$ The proof can easily be modified to show that a homogeneous space
which is a homological Bazaikin space, must be a Berger space or one
of $\Sph^5\times\Sph^9/S^1_{a,b}$. For this one uses the fact that
a homological Bazaikin space has the rational homotopy type of
$\Sph^9\times \CP^2$ (cf. \cite[Lemma 8.2]{BK}) and hence Sullivan's
rational homotopy theory (cf. \cite{GM}) shows that
$\pi_*(X)\otimes \Ra$ is non-zero (and one dimensional) only in
dimensions $2,5,9$.

$b)$ A similar argument gives a short proof that a homogeneous space
which has the rational cohomology of an Eschenburg space must be an
Aloff-Wallach space or $\Sph^3\times \Sph^5/S^1_{a,b}$. The fourth
homotopy group distinguishes the two classes (cp. \cite{E1}). This was
used in \cite{E1} to show that many positively curved Eschenburg
spaces cannot be homotopy equivalent to a homogeneous space. On the
other hand, there are some that are homotopy equivalent but not
homeomorphic to Aloff-Wallach spaces; see \cite{Sh}.

$c)$ In \cite{WZ} it was shown that the homogeneous spaces
$\Sph^5\times\Sph^9/S^1_{a,b}$ satisfy $s=a^3$. For a Bazaikin space
$\B_q$ with $s$ a cube, there hence exist infinitely many homogeneous
spaces with the same cohomology ring. It happens occasionally for
a positively curved Bazaikin space $\B_q$ that $s$ is a cube. The one
with this property and smallest $s$ is given by
$\bar q=(5,5,5,13,23)$, which has $s=17^3$.

$d)$ We finally observe that a Bazaikin space may have the same homotopy
invariants (as above) as the Berger space, e.g.
$\bar{q}=(-89,15,21,31,111)$.
\end{rem*}

\section{Some computer experiments}  

We wrote a C++ program in order to compute all positively curved
Bazaikin spaces $\B_q$ together with the invariants considered in
this paper, and with order $s=|H^6(\B_q)|\leq 40.000.000$. There are
472.959.576 manifolds in a $+20$GB file. Among them, we looked for
subsets of spaces matching different invariants:

\begin{itemize}
\item[$i)$] 3.210.637 pairs, 38.130 triplets, 645 quadruples, 12 5-tuples
and one 6-tuple match the homeomorphism invariants $s$ and $p_1$.
\item[$ii)$] 3.065 pairs match the homotopy invariants $s, lk$ and
$p_1\!\!\mod 24$. The first matching pair (with smallest $s$) is
$\bar q=(-11,13,45,67,77)$ and $(-43,45,49,61,79)$, with
\mbox{$s=254.941$}, $lk=86.294$ and $p_1\!\!\mod 24=7$. No triplet
matches the three invariants.
\item[$iii)$] 365 pairs match $s, p_1$ and $p_2$. No triplets.
The first matching pair is given by $(-13,33,41,105,137)$ and
$(-3,5,77,83,141)$, with $s=999.437$, $p_1=62.271$ and
$p_2=949.280$. It is interesting that all these pairs, but one,
match the first four symmetric functions $\sigma_i(q)$, $1\leq i
\leq 4$. The pair that does not match the first four symmetric
functions (in fact, none of them) is $(-123,149,197,201,525)$ and
$(-19,21,75,437,437)$ with $s=28.864.757$, $p_1=646.383$,
$p_2=25.993.311$.
\item[$iv)$] No pair matches the homeomorphism invariants $s, lk$ and
$p_1$. In the search for such a pair, we went up to
$s\leq 100.000.000$: 2.130.601.485 positively curved Bazaikin
manifolds in a $+100$GB file. Still, no pair matches the homeomorphism
invariants.
\item[$v)$]
In the case of arbitrary (not necessarily positively curved) Bazaikin
manifolds, we found several pairs that match $s$, $lk$, $p_1$ and
$p_2$, especially for small values of $s$, where $lk$ and $p_2$ have
weak influence. There is even a quadruple matching all four
invariants: $(-53,-11,25,33,77)$, $(-53,-23,25,49,69)$,
$(-53, -3,5,49,73)$ and $(-53,-39,41,49,61)$, with $s=1$, $lk=0$,
$p_1=7.807$ and $p_2=0$.
\item[$vi)$] About 7\% of the examples satisfy the (only necessary)
conditions given in \cite{Ba} for the Eschenburg metric on $\B_q$ to
have positive curvature.
\end{itemize}
\medskip

The source code of the programs can be found at
{http://w3.impa.br/$\sim$luis/bazaikin/}

\providecommand{\bysame}{\leavevmode\hbox to3em{\hrulefill}\thinspace}

\end{document}